\documentclass{amsart}

\usepackage{amsmath}
\usepackage{amssymb}
\usepackage{color}
\usepackage{amsthm}
\usepackage{graphicx}
\usepackage{hyperref}

\theoremstyle{plain}

\theoremstyle{definition}

\theoremstyle{plain}

\definecolor{Noir}{rgb}{0,0,0} 
\definecolor{Blanc}{rgb}{1,1,1} 
\definecolor{Gray}{rgb}{0.5,0.5,0.5} 
\definecolor{Rouge}{rgb}{0.8,0.1,0.1} 
\definecolor{DBleu}{RGB}{51,51,178} 
\definecolor{LBleu}{rgb}{0.85,0.85,1} 
\definecolor{Orange}{RGB}{255,140,0}

\newcommand{\bdoc}{\begin{document}} 
\newcommand{\edoc}{\end{document}} 
 
\newcommand{\bcent}{\begin{center}} 
\newcommand{\ecent}{\end{center}} 
 
\newcommand{\benum}{\begin{enumerate}} 
\newcommand{\eenum}{\end{enumerate}} 
\newcommand{\bitem}{\begin{itemize}} 
\newcommand{\eitem}{\end{itemize}} 
 
\newcommand{\btab}{\begin{tabular}} 
\newcommand{\etab}{\end{tabular}} 
\newcommand{\beqn}{\begin{eqnarray}} 
\newcommand{\eeqn}{\end{eqnarray}} 
 
\newcommand{\bmath}{\begin{math}} 
\newcommand{\emath}{\end{math}} 
 
\newcommand{\noin}{\noindent} 
 
\providecommand{\tb}[1]{\textbf{#1}} 
\providecommand{\mb}[1]{\mathbf{#1}} 
 
\newcommand{\bsh}{\backslash} 
 
\newcommand{\ds}{\displaystyle} 
 
\newcommand{\ub}{\underbrace} 
\newcommand{\ob}{\overbrace}

\providecommand{\F}[1]{\mathbb{#1}} 
 
\newcommand{\FF}{\F F} 
\providecommand{\Fn}[1]{\FF_{#1}} 
\newcommand{\Fp}{\Fn{p}} 
\newcommand{\Fq}{\Fn{q}} 
\newcommand{\Fpm}{\Fn{p^m}} 
\newcommand{\Fpn}{\Fn{p^n}} 
\newcommand{\Fpr}{\Fn{p^r}} 
 
\newcommand{\ZZ}{\mathbb{Z}} 
\providecommand{\Zn}[1]{\ZZ_{#1}} 
\providecommand{\ZnZ}[1]{\ZZ/#1\ZZ} 
\newcommand{\Zp}{\Zn{p}} 
\newcommand{\ZpZ}{\ZnZ{p}} 
 
\newcommand{\NN}{\F{N}} 
\newcommand{\QQ}{\F{Q}} 
\newcommand{\RR}{\F{R}} 
\newcommand{\CC}{\F{C}} 
\newcommand{\QQbar}{\overline{\QQ}} 
\newcommand{\Zero}{\mathbb{00}} 
 
\newcommand{\EE}{\F E} 
\newcommand{\II}{\F I} 
\newcommand{\KK}{\F K} 
\newcommand{\MM}{\F M} 
\newcommand{\XX}{\F X} 
\newcommand{\PP}{\F P} 
\newcommand{\FA}{\F A} 
\newcommand{\LL}{\F L} 
\newcommand{\HH}{\mathbb H}  
\newcommand{\FS}{\F S} 
\newcommand{\TT}{\F T}  
\newcommand{\FSig}{\F\Sig} 
\newcommand{\FDel}{\F\Del} 
 
\providecommand{\E}[1]{\hat{\F{#1}}} 
\newcommand{\EC}{\E{C}} 
 
\newcommand{\bQ}{\mathbf{Q}} 
\newcommand{\bP}{\mathbf{P}}

\newcommand{\ind}{\mbox{ind}} 
 
\newcommand{\fx}{f(x)} 
\newcommand{\gx}{g(x)} 
 
\newcommand{\x}{^\star} 
\newcommand{\xs}{^{~\star}} 
\providecommand{\U}[1]{\left(#1\right)\x} 
\newcommand{\xt}{^\times} 
 
\newcommand{\iso}{\simeq} 
\newcommand{\plus}{\oplus} 
\newcommand{\Plus}{\bigoplus} 
\newcommand{\tensor}{\otimes} 
\newcommand{\Tensor}{\bigotimes} 
\newcommand{\inject}{\hookrightarrow} 
\newcommand{\linject}{\hookleftarrow} 
\newcommand{\surject}{\twoheadrightarrow} 
\newcommand{\tri}{\vartriangleleft} 
 
\renewcommand{\ker}{\mbox{ker}\;} 
\newcommand{\Imf}{\mbox{im}~\;} 
\newcommand{\img}{\mbox{im}\;} 
\newcommand{\Hom}{\mbox{Hom}} 
\newcommand{\Sym}{\mbox{Sym}} 
\newcommand{\End}{\mbox{End}\,} 
\newcommand{\Endz}{\mbox{End}_0} 
\newcommand{\Id}{\mbox{Id}} 
\newcommand{\rk}{\mbox{rk}\,} 
\newcommand{\Pic}{\mbox{Pic}} 
\newcommand{\Jac}{\mbox{Jac}} 
\newcommand{\ch}{\mbox{ch}\,} 
\newcommand{\td}{\mbox{td}\,} 
 
\providecommand{\Gal}[1]{\mbox{Gal}(#1)} 
\providecommand{\GAL}[2]{\mbox{Gal}(#1/#2)} 
\providecommand{\Sub}[1]{\mbox{Sub}(#1)} 
\providecommand{\Lat}[1]{\mbox{Lat}(#1)}

\newcommand{\dup}{d_\wedge} 
\newcommand{\drt}{d_>} 

\newcommand{\MF}{\mathfrak}

\newcommand{\Div}{\,\,|\,\,} 
\newcommand{\lcm}{\mbox{lcm}} 
 
\providecommand{\leg}[2]{\left(\frac{#1}{#2}\right)} 
\providecommand{\jac}[2]{\leg{#1}{#2}} 
\providecommand{\qdc}[2]{\left[\frac{#1}{#2}\right]} 
 
\newcommand{\w}{\omega} 
\newcommand{\W}{\Omega} 
\providecommand{\Cal}[1]{\mathcal{#1}} 
\newcommand{\CL}{\Cal L} 
\newcommand{\CO}{\Cal O} 
\renewcommand{\O}{\Cal O} 
\renewcommand{\o}{\Cal O} 
\newcommand{\co}{\Cal O} 
\newcommand{\ca}{\Cal A}
\newcommand{\CR}{\Cal R} 
\newcommand{\CE}{\Cal E} 
\newcommand{\CF}{\Cal F} 
\newcommand{\CQ}{\Cal Q} 
\newcommand{\CK}{\Cal K} 
\newcommand{\CM}{\Cal M} 
\newcommand{\CN}{\Cal N} 
\newcommand{\CDee}{\Cal D} 
\newcommand{\CP}{\Cal P} 
\newcommand{\CS}{\Cal S} 
\newcommand{\ClC}{\Cal C} 
\newcommand{\CJ}{\Cal{J}} 
\newcommand{\CA}{\Cal{A}} 
\newcommand{\CB}{\Cal{B}}

\newcommand{\Si}{\Sigma} 
 
\providecommand{\Ok}[1]{\CO_{#1}} 
\newcommand{\OK}{\Ok{K}}

\renewcommand{\epsilon}{\varepsilon} 
\newcommand{\ep}{\varepsilon} 
 
\providecommand{\abs}[1]{\left|#1\right|} 
\providecommand{\norm}[1]{\lVert#1\rVert} 
 
\newcommand{\di}{\partial} 
\providecommand{\ddy}[1]{\ds\frac{d}{d #1}} 
\providecommand{\didiy}[1]{\ds\frac{\di}{\di #1}} 
\newcommand{\ddx}{\ddy{x}} 
\newcommand{\didix}{\didiy{x}} 
\providecommand{\v}[1]{\vec{#1}} 
 
\newcommand{\ihat}{\hat{\infty}}

\providecommand{\set}[1]{\left\{#1\right\}} 
\providecommand{\lst}[1]{\ds\left[\,\,#1\,\,\right]}

\providecommand{\ip}[2]{\left( #1,#2\right)} 
\providecommand{\IP}[2]{\left\langle #1,#2\right\rangle} 
\providecommand{\bra}[1]{\left\langle\left. #1\right|\right.} 
\providecommand{\ket}[1]{\left.\left| #1\right.\right\rangle} 
\providecommand{\bkv}[4]{\left\langle\left.\tb{#1} #2\right|\tb{#3} #4\right\rangle} 
\providecommand{\bk}[2]{\bkv{}{#1}{}{#2}} 
\providecommand{\bkvop}[5]{\left\langle\tb{#1} #2\left| #5\right|\tb{#3} #4\right\rangle} 
\providecommand{\bkop}[3]{\bkvop{}{#1}{}{#2}{#3}} 
\providecommand{\comm}[2]{\left[#1,#2\right]} 
\providecommand{\expn}[1]{\left\langle #1\right\rangle} 
\providecommand{\gen}[1]{\expn{#1}} 
\newcommand{\Del}{\Delta} 
\newcommand{\Nab}{\nabla} 
\newcommand{\Sig}{\Sigma} 
\newcommand{\oline}{\overline}

\newcommand{\p}{\rho} 
\providecommand{\sprod}[2]{\left\langle #1,#2\right\rangle}

\newcommand{\Ehat}{\overline E} 
\newcommand{\Phihat}{\overline \Phi} 
\newcommand{\phihat}{\overline \phi}

\newcommand\varleq{\mathbin{\vcenter{\hbox{%
  \oalign{\hfil$\scriptstyle<$\hfil\cr 
          \noalign{\kern-.3ex} 
          $\scriptscriptstyle({-})$\cr}%
}}}} 
 
\renewcommand\subsetneq{\mathbin{\vcenter{\hbox{%
  \oalign{\hfil$\scriptstyle\subset$\hfil\cr 
          \noalign{\kern-.3ex} 
          $\scriptscriptstyle({-})$\cr}%
}}}} 

\author{Aidan Lindberg}
\address{Department of Mathematics \& Statistics\\
  McLean Hall, University of Saskatchewan\\
  Saskatoon, SK, Canada  ~S7N 5E6\vspace{10pt}}
\email{aidan.lindberg@usask.ca, rayan@math.usask.ca}

\author{Steven Rayan}

\title[Geodesics on a KN(A)dS instanton]{Geodesics on a Kerr-Newman-(anti-)de Sitter instanton}
\date{\today}

\begin{document}

\maketitle

\noin\tb{Abstract.} We study geodesics along a noncompact Kerr-Newman instanton, where the asymptotic geometry is either de Sitter or anti-de Sitter.  We use first integrals for the Hamilton-Jacobi equation to characterize trajectories both near and away from horizons.  We study the interaction of geodesics with special features of the metric, particularly regions of angular degeneracy or ``theta horizons'' in the de Sitter case.  Finally, we characterize a number of stable equilibrium orbits.\\

\tableofcontents

\section{Introduction}

Much is known already regarding geodesics along Kerr spacetimes.  Physically, these correspond to particles under the influence of a rotating black hole.  Such geodesics have also been studied in the presence of charge and a nonzero cosmological constant (e.g. \cite{Rayan}).  Less studied are geodesics along the corresponding gravitational instanton.  In the Riemannian setting, there are no timelike curves that play the role of usual matter.  That being said, geodesics in gravitational instantons do have physical interpretations: for instance, as the relative motion of certain monopoles.  Notably, the Kaluza-Klein monopole can be embedded in a $5$D Taub-NUT instanton, \cite{Visi}.  We will treat only a $4$D instanton (of Kerr type) here, but the differential-geometric question of classifying its geodesics is still well-posed and may have applications to various scattering questions (monopole or otherwise) in instanton backgrounds.  Generally speaking, Riemannian solutions to the Einstein equations have emerged as state transition probabilities in Euclidean quantum gravity \cite{BoothMann,Chrusciel,Chrusciel2} and are expected to encode quantum properties of their Lorentzian counterparts.  Combining this with the general importance of the Kerr solution in the context of recent developments such as AdS/CFT and Kerr/CFT, understanding the geometric structure of the Riemannian analogue of Kerr is a priority.  

Here, we consider a noncompact, incomplete Kerr-Newman-(anti-)de Sitter instanton.  One can employ a periodic identification of the imaginary time coordinate to remove the singularity and create a compact version of the instanton, as suggested in \cite{GibbonsHawking,GibbonsHawking2} and then carried out and explored in \cite{Chrusciel,Chrusciel2}.  This compact space will have a complete, everywhere positive-definite Einstein metric.   We, on the other hand, obtain the instanton metric from a Wick rotation but do not perform the periodic identification.  As such, the solution has degeneracies and signature changes.  In particular, it is only Riemannian on certain submanifolds, and so one caveat is that the solution is a gravitational instanton only in a weaker sense.

Given that geodesics along the manifold in question possess fixed values of the first integrals that correspond usually to rest mass, energy, and angular momentum, we refer to them in any event as ``particles'' (even if they originate in the Riemannian part of the manifold, where they are spacelike).    For both the de Sitter and anti-de Sitter geometries, we study the effect of the first integrals on the geodesic evolution, paying particular attention to their potential functions and to their sensitivity to the cosmological constant, the charge, and whether or not the singularity is slowly or rapidly rotating.  In the de Sitter situation, the metric degenerates in an interesting fashion --- there are surfaces of angular degeneracy (as referred to in \cite{Stuk}) --- and we study the trajectories of these particles near these degenerations.  We also characterize a number of stable equilibrium orbits in this case.

\subsection{Metric}

Here, we consider the Kerr-Newman-(anti)-de Sitter, or KN(A)dS, instanton, which has underlying $4$-manifold $M=\RR^2\times S^2$ and a metric $g$ that solves the Riemannian Einstein-Maxwell field equations with positive (negative) cosmological constant.  The field equations are\begin{eqnarray}R_{\mu\nu}+\left(\Lambda-\frac{1}{2}R\right)g_{\mu\nu} & = & 2\left(F_{\mu\alpha}F^\alpha_\nu-\frac{1}{4}F^{\alpha\beta}F_{\alpha\beta}g_{\mu\nu}\right)\nonumber\\dF\;=\;d\star F & = & 0\nonumber,\end{eqnarray} where $g_{\mu\nu}$ is the metric, $R_{\mu\nu}$ and $R$ are the Ricci tensor and scalar curvature respectively, $\Lambda$ is the cosmological constant (a fixed real number), $F$ is the Maxwell field, and $\star$ is the Hodge star.  Choosing $\Lambda$ positive results in a geometry that is asymptotically de Sitter while $\Lambda$ negative results in an anti-de Sitter geometry.

After fixing $\Lambda$ (and setting the magnetic charge parameter $p$ to $0$), there is a family of Lorentzian KN(A)dS metrics depending on $3$ real numbers: $M$, the total mass of the instanton; $a$, the angular momentum per unit mass; and $e$, the charge per unit mass.  That such solutions to the Einstein-Maxwell equations are parametrized by $a,e,M$ is essentially the ``No Hair Theorem'' of general relativity.  The KN(A)dS instanton solution is obtained from the corresponding spacetime solution in the standard way via a Wick transformation:$$e\mapsto ie,\;a\mapsto ia,\mbox{ and }t\mapsto it,$$ where $i=\sqrt{-1}$ (cf. \cite{Chrusciel,Chrusciel2} for instance).  The Wick-rotated metric $g$ then has the following line element in Boyer-Lindquist coordinates $(r,t,\theta,\varphi)$ on $\RR^2\times S^2$:
\begin{equation}
ds^2=\frac{\Sigma}{\Delta_r}dr^2+\frac{\Sigma}{\Delta_\theta}d\theta^2+\frac{S^2}{\Xi^2\Sigma}\Delta_\theta(a dt+(r^2-a^2)d\varphi)^2+\frac{\Delta_r}{\Xi^2\Sigma}(dt-aS^2 d\varphi)^2,\nonumber
\end{equation}
in which we have functions
\begin{align}
\Sigma=&r^2-a^2\cos^2\theta\nonumber\\
\Delta_r=&(r^2-a^2)(1-Lr^2)-2Mr-e^2\nonumber\\
\Delta_\theta=&1-La^2\cos^2\theta\nonumber\\
\Xi=&1-La^2\nonumber, 
\end{align}
with $L=\displaystyle\frac{\Lambda}{3}$.  We also note that the Maxwell potential under which $g$ solves the field equations is
\begin{equation}
A=\frac{er}{\Sigma \Xi}(-dt+a\sin^2\theta d\varphi).\nonumber
\end{equation}

We partition the instanton in accordance with the roots of the $\Delta_r$ function. There are exactly two real roots for the anti-de Sitter instanton and four for its de Sitter counterpart (except for at extreme values of the parameters, as discussed in the Appendix). The partitions are defined in Table 1.1 and Table 1.2.\\

\begin{center}
Table 1.1: \textit{Block Structure of AdS $\Delta_r$ Roots}\vspace{5mm}  \linebreak
\begin{tabular}{c|c}
$r\in (r_+,\infty)$ & AdS Block \\ \hline
$r\in (r_-,r_+)$ & Block I \\ \hline
$r\in (-\infty, r_-)$ & Block II \\
\end{tabular}
\end{center}
\vspace{5pt}
\begin{center}
Table 1.2: \textit{Block Structure of dS $\Delta_r$ Roots}\vspace{5mm}  \linebreak
\begin{tabular}{c|c}
$r\in (r_{++},\infty)$ & dS$_+$ Block \\ \hline
$r\in (r_+,r_{++})$ & Block I \\ \hline
$r\in (r_- r_+)$ & Block II \\ \hline
$r\in (r_{--},r_-)$ & Block III \\ \hline
$r\in (-\infty,r_{--})$ & dS$_-$ Block
\end{tabular}
\end{center}
\vspace{5pt}
Note that there will only exist one negative root for both instantons. This is $r_-$ in the AdS instanton, and $r_{--}$ in the dS instanton. The details of determining these root structures are deferred to the Appendix.

It is worth noting that, while we use the word ``instanton'' to describe the metric, this is only a gravitational instanton in an incomplete sense: we do not impose the ALE decay at infinity and the metric signature is not constant.  (This is in constrast to \cite{Chrusciel,Chrusciel2}, who construct a compact, singularity-free version of the solution in which there is constant Riemannian signature.)  This can be seen in the eigenvalues of the metric: 
\begin{align*}
\lambda_1&=g_{\theta \theta}\nonumber \\
\lambda_2&=g_{r r}\nonumber \\
\lambda_3&=\frac{g_{\varphi \varphi}+g_{t t}+\sqrt{g_{\varphi \varphi}^2-2g_{\varphi \varphi}g_{t t}+g_{tt}^2+4g_{t \varphi}^2}}{2}\nonumber \\
\lambda_4&=\frac{g_{\varphi \varphi}+g_{t t}-\sqrt{g_{\varphi \varphi}^2-2g_{\varphi \varphi}g_{t t}+g_{tt}^2+4g_{t \varphi}^2}}{2}\nonumber.
\end{align*}
The positivity of these functions requires the following inequality to be satisfied:
\begin{align*}
0&<\frac{(r^2-a^2)^2}{a^2S^2}+a^2S^2+2(r^2-a^2).
\end{align*}

Thus, for the anti-de Sitter instanton, we have the following signatures:
\begin{itemize}
\item \textbf{Riemannian}: Outside of the singularity, and in Block II or the AdS Block
\item \textbf{Lorentzian}: Inside of Block I
\item \textbf{Two Positive and Two Negative Eigenvalues}: Inside of the singularity and in Block II or the AdS Block
\end{itemize}
and for the de Sitter instanton, we have the following signatures:
\begin{itemize}
\item \textbf{Riemannian}:
	\begin{itemize}
	\item Inside of the singularity: Inside of the cones and in Block II or in dS$_\pm$
	\item Outside of the singularity: Outside of the cones and in Block I or Block III
	\end{itemize}
\item \textbf{Lorentzian}:
	\begin{itemize}
	\item Outside of the cones and in Block II or dS$_\pm$
	\item Inside of the cones and in Block I or Block III
	\end{itemize}
\item \textbf{Two Positive and Two Negative Eigenvalues}:
	\begin{itemize}
	\item Inside of the singularity: Outside of the cones and in Block I or Block III
	\item Outside of the singularity: Inside of the cones and in Block II or dS$_\pm$
	\end{itemize}
\end{itemize}

\subsection{Canonical Vector Fields}

As is clear from the line element, the vector fields $\di_t$ and $\di_\phi$ are not orthogonal.  We can correct this by introducing new vector fields $V$ and $W$ (as in \cite{ONeill}, for instance):
\begin{align*}
V:=&(r^2-a^2)\partial_t-a\partial_\varphi \\
W:=&\partial_\varphi+aS^2\partial_t.
\end{align*}
We now have an orthogonal basis for the space of vector fields on $M$, namely$$\{ V, W, \partial_r, \partial_\theta\}.$$ These yield the following standard identities, which will be useful when deriving the equations of motion for test particles:
\begin{align*}
\langle V, \partial_\varphi \rangle =(r^2-a^2)g_{\varphi t}-a^2g_{\varphi \varphi}&=\frac{-aS^2\Delta_r}{\Xi^2} \\
\langle V, \partial_t \rangle = (r^2-a^2)g_{t t}-ag_{t\varphi}&=\frac{\Delta_r}{\Xi^2} \\
\langle W, \partial_\varphi \rangle = g_{\varphi \varphi}+aS^2 g_{t \varphi}&=\frac{S^2(r^2-a^2)\Delta_\theta}{\Xi^2} \\
\langle W, \partial_t \rangle=g_{t \varphi}+aS^2g_{tt}&=\frac{aS^2\Delta_\theta}{\Xi^2} \\
\langle V, V \rangle &= \frac{\Delta_r \Sigma}{\Xi^2} \\
\langle W, W \rangle &= \frac{S^2\Delta_\theta \Sigma}{\Xi^2}.
\end{align*}

\subsection{Singularity, $r$-horizons, and $\theta$-horizons}

In the AdS instanton, the function $\Delta_r$ will always have two real roots, as well as two complex (non-physical) roots. Also, it is worth noting that there is no possibility of $\theta$-horizons (called ``surfaces of degeneracy'' in \cite{Stuk}) for $\Lambda<0$. For the dS instanton, we will have four real roots if we impose reasonable assumptions on parameters (see Appendix), and $\theta$-horizons provided $a>\sqrt{\frac{3}{\Lambda}}$. Perhaps the most interesting difference between this metric and the ordinary KN(A)dS metric is its peculiar singularity structure. Our $\Sigma$ function will be zero whenever
\begin{equation*}
\frac{r}{a}=\cos\theta.
\end{equation*}
This induces a singularity which takes on the radial values $r\in[-a,a]$. It must collapse to a point at both $r=0$ and $r=\pm$a, however, it may take on multiple colatitude values when $r\in(-a,0)\cup(0,a)$, creating a geometry resembling two spheres whose poles have been stretched out and joined at $r=0$. This singularity will be present within the Riemannian submanifold provided $|r_-|\leq a$ or $r_+\leq a$.

\section{Integrability of the Equations of Motion}

The natural Lagrangian for test particles (e.g. \cite{Sharp}) is:
\begin{equation}
\textbf{L}_q(s,x^j,\dot{x}^j)=\frac{1}{2}\textit{g}_{ab}\dot{x}^a\dot{x}^b-qA_a\dot{x}^a.
\end{equation}
Instead of working directly with this second-order equation, we go down the well-trodden road of studying sufficiently-many \emph{first integrals} in involution with the corresponding Hamiltonian, which we denote by $\textbf{H}$, as per standard integrability theory in symplectic geometry (cf. \cite{Lou}).  These can be found by separating variables appropriately for the Hamilton-Jacobi equation.  For various Kerr-type backgrounds, this is well-known, with the result appearing originally in \cite{Carter}.  A general result of this sort can be established for a large class of field-equation solutions via \cite{Kamran}.  For the sake of completeness, and to fix the form of the appropriate first-order equations of motion, we derive the separation for the Kerr-Newman-(anti)-de Sitter background directly here, regardless of signature.   We adapt the method of \cite{ONeill}, where the separation is demonstrated for the pure, Lorentzian Kerr solution (i.e. $e=\Lambda=0$).  (For nonzero $e$ and positive $\Lambda$, an adaptation of this same technique is used in \cite{Rayan}.)   We find the existence of the two Killing vector fields $\partial_\varphi$ and $\partial_t$, along with the Hamiltonian itself, giving us three of the four necessary integrals:
\begin{itemize}
\item $\mathfrak{q}:=2\textbf{H}$ is the particle's \textit{rest mass}\\
\item $\mathcal{E}:=-\textit{p}_t$ is the \textit{energy} of the particle\\
\item $\mathcal{L}:=\textit{p}_\varphi$ is the \textit{angular momentum} of the particle\\
\end{itemize}
The fourth first integral is typcally referred to as \emph{Carter's constant}.  To isolate it here, we follow \cite{ONeill} and examine the quantity
\begin{equation*}
\Psi=\frac{\langle \gamma',V \rangle}{\langle V, V \rangle}V+\frac{\langle \gamma', W \rangle}{\langle W, W \rangle}W,
\end{equation*}
which can be expanded as 
\begin{align*}
\Psi&=\frac{t'\langle \partial_t, V \rangle +\varphi' \langle \partial_\varphi, V \rangle}{\Delta_r/\Xi^2}V+\frac{t' \langle \partial_t, W \rangle +\varphi' \langle \partial_\varphi, W \rangle}{S^2\Delta_\theta/\Xi^2}W \\
\Psi&=\frac{t'(V+aW)+\varphi'((r^2-a^2)W-aS^2V)}{\Sigma} \\
\Psi&=t'\partial_t+\varphi'\partial_\varphi.
\end{align*}
Thus, we may write
\begin{equation}
\gamma'=r'\partial_r+\theta'\partial_\theta+\frac{\langle \gamma\, V \rangle}{\langle V, V \rangle}V+\frac{\langle \gamma',W \rangle}{\langle W, W \rangle}W,
\end{equation}
and
\begin{align*}
\langle \gamma', \gamma' \rangle&=r'^2\frac{\Sigma}{\Delta_r}+\theta'^2\frac{\Sigma}{\Delta_\theta}+\frac{\langle \gamma',V\rangle^2}{\Delta_r/\Xi^2}+\frac{\langle \gamma',W\rangle^2}{S^2\Delta_\theta/\Xi^2} \\
\mathfrak{q}\Sigma&=\frac{r'^2\Sigma^2+\Xi^2\langle\gamma',V\rangle^2}{\Delta_r}+\frac{\theta'^2\Sigma^2+\Xi^2\langle\gamma',W\rangle^2/S^2}{\Delta_\theta}.
\end{align*}
Note that the left hand side of this equation is well defined everywhere in the instanton, thus the right hand side must be as well. This can be expanded and rearranged to give
\begin{equation*}
\mathcal{K}:=qr^2-\frac{r'^2\Sigma^2+\Xi^2\langle\gamma',V\rangle^2}{\Delta_r}=qa^2C^2+\frac{\theta'^2\Sigma^2+\Xi^2\langle\gamma',W\rangle^2/S^2}{\Delta_\theta}.
\end{equation*}
Writing the Hamiltonian with the contravariant metric in the form:
\begin{equation}
\textbf{H}=\frac{1}{2}g^{ab}(p_a+qA_a)+(p_b+qA_b),
\end{equation}
we suppose that the Hamiltonian in this form may be separated into
\begin{equation*}
\textbf{H}=\frac{1}{2}\frac{H_r+H_\theta}{U_r+U_\theta},
\end{equation*}
where $U_r$ and $U_\theta$ are single variable functions of \textit{r} and $\theta$ respectively, and we have:
\begin{itemize}
\item $H_r$ being independent of $p_\theta$ and all Boyer-Lindquist coordinates, save for r
\item $H_\theta$ being independent of $p_r$ and all Boyer-Lindquist coordinates, save for $\theta$
\end{itemize}
It is proven in \cite{Carter} that, should the Hamiltonian separate into the above form, the function 
\begin{equation*}
\mathcal{K}=\frac{U_r H_\theta-U_\theta H_r}{U_r+U_\theta}
\end{equation*} 
is in involution with \textbf{H} and is a first integral. 

To show the desired separation, we first examine the contravariant metric 
\[ g^{ab}=
\begin{bmatrix}
\frac{1}{\alpha} & 0 & 0 & 0 \\
0 & \frac{1}{\beta} & 0 & 0 \\
0 & 0 & \frac{\psi}{d} & \frac{-\sigma}{d} \\
0 & 0 & \frac{-\sigma}{d} & \frac{\eta}{d}
\end{bmatrix}
\]
where we have $\eta=g_{tt}, \sigma=g_{t \varphi}, \psi=g_{\varphi \varphi}$ and
\begin{equation*}
d=\psi\eta-\sigma^2=\frac{S^2\Delta_\theta \Delta_r}{\Xi^4}.
\end{equation*}
This determinant is non-zero everywhere except for the horizons and poles, where Boyer-Lindquist coordinates fail in any case. This allows us to rewrite the Hamiltonian as
\begin{equation*}
H=\frac{1}{2}\left[\frac{p_r^2}{g_{rr}}+\frac{p_{\theta}^2}{g_{\theta \theta}}+\frac{g_{\varphi \varphi}(p_t+qA_t)^2}{d}+\frac{g_{tt}(p_\varphi+qA_\varphi)^2}{d}\right]-\frac{g_{t \varphi}(p_t+qA_t)(p_\varphi+qA_\varphi)}{d}.
\end{equation*}
Inserting the metric components and simplifying yields:
\begin{multline*}
2\Sigma H=p_r^2\Delta_r+p_\theta^2\Delta_\theta+\Xi^2\left(\frac{(r^2-a^2)^2}{\Delta_r}+\frac{a^2S^2}{\Delta_\theta}\right)\left(p_t-\frac{qer}{\Sigma\Xi}\right)^2+ \\\Xi^2\left(\frac{a^2}{\Delta_r}+\frac{1}{S^2\Delta_\theta}\right)\left(\frac{qeraS^2}{\Sigma\Xi}+p_\varphi\right)^2-\\ 2a\Xi^2\left(\frac{r^2-a^2}{\Delta_r}-\frac{1}{\Delta_\theta}\right)\left(p_t-\frac{qer}{\Sigma\Xi}\right)\left(\frac{qeraS^2}{\Sigma\Xi}+p_\varphi\right).
\end{multline*}
Further simplification yields the desired separated functions:
\begin{align*}
H_r(r,p_r,p_t,p_\varphi)&=p_r^2\Delta_r+\frac{\Xi^2}{\Delta_r}\left(\frac{qer}{\Xi}+ap_\varphi-(r^2-a^2)p_t\right)^2 \\
H_\theta(\theta, p_\theta, p_t, p_\varphi)&=p_\theta^2\Delta_\theta+\frac{\Xi^2}{S^2\Delta_\theta}\left(aS^2p_t+p_\varphi\right)^2 \\
U_r(r)&=r^2 \\
U_\theta(\theta)&=-a^2C^2.
\end{align*}
Evidently $U_r+U_\theta=\Sigma$ and $H_r+H_\theta=2\Sigma H$ as required. We then examine the function
\begin{align*}
\mathcal{K}&=\frac{U_rH_\theta-U_\theta H_r}{U_r+U_\theta} \\
&=\frac{U_rH_\theta-U_\theta H_r+(U_rH_r-U_rH_r)}{U_r+U_\theta} \\
&=\mathfrak{q}U_r-H_r.
\end{align*}
Then, with the definition
\begin{equation*}
\mathcal{P}:=\frac{qer}{\Xi}+ap_\varphi-(r^2-a^2)p_t,
\end{equation*}
we can see
\begin{equation*}
\mathcal{K}=qr^2-\Delta_rp_r^2-\frac{\Xi^2\mathcal{P}^2}{\Delta_r}.
\end{equation*}
At this point, we transfer to tangent space coordinates by $x_k'=\frac{\partial H}{\partial p'^k}$, thus replacing $p_r$ by $g_{rr}r'-qA_r=\frac{\Sigma}{\Delta_r}r'$. Substituting into $\mathcal{K}$ gives
\begin{equation}
\mathcal{K}=\mathfrak{q}r^2-\frac{\Sigma^2r'^2+\Xi^2\mathcal{P}^2}{\Delta_r}.
\end{equation}
Next, we examine $\mathcal{P}(r)$ in tangent space coordinates
\begin{align*}
\mathcal{P}&=\frac{qer}{\Xi}+a(g_{\varphi \varphi} \varphi'+g_{t \varphi}t'-qA_\varphi)-(r^2-a^2)(g_{tt}t'+g_{t \varphi}\varphi'-qA_t) \\
&=\langle-(r^2-a^2)\partial_t, t'\partial_t+\varphi'\partial_\varphi\rangle+\langle a \partial_\varphi, \varphi'\partial_\varphi+t'\partial_t\rangle+q\left(\frac{er}{\Xi}-aA_\varphi+(r^2-a^2)A_t\right) \\
&=-\langle t'\partial_t+\varphi' \partial_\varphi, V\rangle+q\left(\frac{er}{\Xi}-\frac{a^2S^2er}{\Sigma\Xi}+\frac{era^2}{\Sigma\Xi}-\frac{er^2}{\Sigma\Xi}\right) \\
&=-\langle t'\partial_t+\varphi' \partial_\varphi, V\rangle.
\end{align*}
Using our orthogonal basis, we find
\begin{equation*}
\mathcal{P}=-\langle r'\partial_r+\theta'\partial_\theta+\varphi'\partial_\varphi+t'\partial_t,V\rangle=-\langle \gamma',V\rangle.
\end{equation*}
This give us 
\begin{equation}
\mathcal{K}=\mathfrak{q}r-\frac{\Sigma^2r'^2+\Xi^2\langle \gamma',V\rangle^2}{\Delta_r},
\end{equation}
as required. Alternatively, we may use our $H_\theta$ function and define $\mathcal{D}(\theta):=(a^2S^2p_t+p_\varphi)=\langle \gamma',W\rangle$. Taking $\mathcal{K}=H_\theta-\mathfrak{q}U_\theta$, and again transforming to tangent space coordinates, we see that 
\begin{equation*}
\mathcal{K}=p_\theta^2\Delta_\theta+\frac{\Xi^2\mathcal{D}^2}{\Delta_\theta S^2},
\end{equation*}
or equivalently, 
\begin{equation}
\mathcal{K}=\mathfrak{q}a^2C^2+\frac{\Sigma^2\theta'^2+\Xi^2\langle\gamma',W\rangle^2/S^2}{\Delta_\theta},
\end{equation}
giving us the alternate form of Carter's constant, and showing that this is indeed our fourth first integral. This gives the complete set of first integrals for the instanton $\{\mathfrak{q}, \mathcal{E}, \mathcal{L}, \mathcal{K}\}$.
Using the functions $\mathcal{P}$ and $\mathcal{D}$ to represent $\mathcal{K}$, a slight rearrangement yields the two first integral equations of motion:
\begin{align}
R(r):=&\Sigma^2r'^2=\Delta_r(\mathfrak{q}r^2-\mathcal{K})-\Xi^2\mathcal{P}^2 \\
\Theta(\theta):=&\Sigma^2\theta'^2=\Delta_\theta(\mathcal{K}-\mathfrak{q}a^2C^2)-\frac{\Xi^2\mathcal{D}^2}{S^2}.
\end{align}
These are our \textit{radial} and \textit{colatitude equations} respectively. To find our last two first integral equations, we start with the first order coupled equations:
\begin{align*}
-\mathcal{E}+\mathfrak{q}A_t&=t'g_{tt}+\varphi'g_{t\varphi} \\
\mathcal{L}+\mathfrak{q}A_\varphi &=t'g_{t\varphi}+\varphi'g_{\varphi \varphi}.
\end{align*}
Substituting multiples of these equations by metric terms, so as to remove $\varphi'$ dependence gives us 
\begin{equation*}
t'(g_{tt}g_{\varphi\varphi}-g_{t\varphi}^2)=-g_{\varphi\varphi}\mathcal{E}-g_{t\varphi}\mathcal{L}+q(g_{\varphi\varphi}A_t-g_{t\varphi}A_\varphi).
\end{equation*}
Noting that the quantity on the LHS is the determinant of the contravariant metric, and substituting in terms, we find
\begin{equation*}
\frac{t'\Delta_\theta\Delta_r\Sigma}{\Xi^2}=\Delta_\theta (r^2-a^2)\left(-\mathcal{E}(r^2-a^2)-\mathcal{L}a-\frac{qer}{\Xi}\right)+\Delta_r a(\mathcal{L}-a\mathcal{E}S^2).
\end{equation*}
Solving for $t'$ gives
\begin{equation*}
t'=\frac{\Xi^2}{\Sigma}\left(\frac{-(r^2-a^2)\mathcal{P}}{\Delta_r}+\frac{\mathcal{D}a}{\Delta_\theta}\right).
\end{equation*}
Repeating this procedure, but substituting multiples of metric terms to remove $t'$ dependence yields
\begin{equation*}
\frac{\varphi' S^2 \Delta_\theta \Delta_r \Sigma}{\Xi^2}=aS^2\Delta_\theta\left((r^2-a^2)\mathcal{E}+a\mathcal{L}+\frac{qer}{\Xi}\right)+\Delta_r(\mathcal{L}-a\mathcal{E}S^2).
\end{equation*}
Solving for $\varphi'$ then shows
\begin{equation*}
\varphi'=\frac{\Xi^2}{\Sigma}\left(\frac{a\mathcal{P}}{\Delta_r}+\frac{\mathcal{D}}{S^2\Delta_\theta}\right).
\end{equation*}
Thus, we have our final two equations:
\begin{align}
t'&=\frac{\Xi^2}{\Sigma}\left(\frac{-(r^2-a^2)\mathcal{P}}{\Delta_r}+\frac{\mathcal{D}a}{\Delta_\theta}\right) \\
\varphi'&=\frac{\Xi^2}{\Sigma}\left(\frac{a\mathcal{P}}{\Delta_r}+\frac{\mathcal{D}}{S^2\Delta_\theta}\right),
\end{align}
referred to as the $t$ and $\varphi$ \textit{first integral equations} respectively. 

\section{KNAdS Instanton Particle Dynamics}
In the Kerr-Newman-anti-de Sitter instanton, there are many similarities to the Kerr-Newman-de Sitter spacetime, owing to the similar root structure of $\Delta_r$ and the absence of $\Delta_\theta=0$ degeneracies. We begin by defining $\gamma'_{\Pi}:=\mbox{Span}\{\partial_r,V\}$ and $\gamma'_\perp:=\mbox{Span}\{\partial_\theta, W\}$. In the case where $\gamma'_\perp=0$ everywhere, we call this a \textit{principal orbit}, as in \cite{ONeill}. Writing out $\gamma'$ in the form
\begin{equation*}
\gamma'=r'\partial_r+\theta' \partial_\theta+\frac{\langle \gamma', W \rangle}{\langle W, W \rangle}W+\frac{\langle \gamma', V \rangle}{\langle V, V \rangle}V,
\end{equation*}
we take product of the components in $\gamma'_\perp$
\begin{equation*}
\langle \gamma'_\perp, \gamma'_\perp \rangle =\theta'^2g_{\theta\theta}+\frac{\mathcal{D}^2}{\langle W,W \rangle}.
\end{equation*}
Solving for $\theta'^2$ and substituting in terms we see
\begin{equation*}
\theta'^2=\frac{\Delta_\theta}{\Sigma}\left(\langle \gamma'_\perp,\gamma'_\perp \rangle-\frac{\mathcal{D}^2\Xi^2}{\Sigma S^2 \Delta_\theta}\right).
\end{equation*}
Substituting this into our expression for $\mathcal{K}$, we find 
\begin{equation*}
\mathcal{K}=\mathfrak{q}a^2C^2+\Sigma\langle\gamma'_\perp,\gamma'_\perp \rangle.
\end{equation*}
Noting that $\mathfrak{q}=\langle\gamma'_\Pi,\gamma'_\Pi\rangle+\langle\gamma'_\perp,\gamma'_\perp\rangle$, we see 
\begin{align*}
\mathcal{K}&=\Sigma\langle\gamma'_\perp,\gamma'_\perp\rangle+\left(\langle\gamma'_\Pi,\gamma'_\Pi\rangle+\langle\gamma'_\perp,\gamma'_\perp\rangle\right)a^2C^2 \\
&=r^2\langle\gamma'_\perp,\gamma'_\perp\rangle+(r^2\langle\gamma'_\Pi,\gamma'_\Pi\rangle-r^2\langle\gamma'_\Pi,\gamma'_\Pi\rangle)+\langle\gamma'_\Pi,\gamma'_\Pi\rangle a^2C^2, \\
\end{align*}
which simplifies to 
\begin{equation}
\mathcal{K}=\mathfrak{q}r^2-\Sigma\langle\gamma'_\Pi, \gamma'_\Pi\rangle.
\end{equation}
One important distinction between the instanton and non-instanton spacetimes is the existence of negative $\Sigma$ terms inside of the singularity, which we denote by $\Omega$. With this in mind, we can study the causal character of $\langle\gamma'_\Pi,\gamma'_\Pi\rangle$. As $\langle V,V\rangle$ and $g_{rr}$ will always have the same sign, we summarize the causal character of $\gamma'_\Pi$ in the following table:

\vspace{10pt}

\begin{center}
Table 3.1: \textit{Causal Character of $\langle\gamma'_\Pi,\gamma'_\Pi\rangle$ in KNAdSI}\vspace{3mm}
\begin{tabular}{c|c}
$\langle\gamma'_\Pi,\gamma'_\Pi\rangle<0$ & AdS Block inside $\Omega$ -- Block I outside $\Omega$ -- Block II inside $\Omega$ \\ \hline
$\langle\gamma'_\Pi,\gamma'_\Pi\rangle>0$ & AdS Block outside $\Omega$ -- Block I inside $\Omega$ -- Block II outside $\Omega$ 
\end{tabular}
\end{center}

\vspace{10pt}

In this $\theta$-horizon-free instanton, the causal character of $\langle\gamma'_\perp,\gamma'_\perp\rangle$ will be timelike inside of $\Omega$ and spacelike outside of $\Omega$. Using the $\gamma'_\perp$ formulation of $\mathcal{K}$, we summarize the possible $\mathcal{K}$ values in the following table:

\vspace{10pt}

\begin{center}
Table 3.2: \textit{Possible $\mathcal{K}$ values for various $\mathfrak{q}$}\vspace{3mm}
\begin{tabular}{c|c}
$\mathfrak{q}<0$ & $\mathcal{K}\geq \mathfrak{q}a$ with equality if a principal orbit at a pole \\ \hline
$\mathfrak{q}=0$ & $\mathcal{K}=0$ if principal and $\mathcal{K}=\Sigma\langle\gamma'_\Pi,\gamma'_\Pi\rangle$ otherwise \\ \hline
$\mathfrak{q}>0$ & $\mathcal{K}\geq 0$ with equality if a principal equatorial orbit
\end{tabular}
\end{center}

\vspace{10pt}

\subsection{Time and Longitude Coordinate Evolution}
We consider the case of particles with $\mathcal{E}=\mathcal{L}=0$ i.e. "\textit{lazy particles}". These lazy particles have their \textit{t} and $\varphi$ coordinate functions simplify to 
\begin{align}
t'&=\frac{-(r^2-a^2)qer\Xi}{\Sigma\Delta_r} \\
\varphi'&=\frac{\Xi aqer}{\Sigma\Delta_r},
\end{align}
respectively. From this, we can construct the following tables (3.3 and 3.4) summarizing how each of these functions will evolve in the various areas of the instanton.

\vspace{10pt}

\begin{center}
Table 3.3: \textit{Direction of travel along the $\varphi'$ axis for $\mathcal{L}=\mathcal{E}=0$ with $\Lambda<0$}\vspace{3mm}
\begin{tabular}{c|c|c|c|c}
\, &AdS&Block I with $r>0$ &Block I with $r<0$ & Block II\\ \hline
$qe<0$ & - & + & - & + \\ \hline
$qe>0$ & + & - & + & - \\ \hline
\end{tabular}
\end{center}

\vspace{10pt}

\begin{center}
Table 3.4: \textit{Direction of travel along the $t'$ axis for $\mathcal{L}=\mathcal{E}=0$ with $\Lambda<0$}\vspace{3mm}
\begin{tabular}{c|c|c|c|c}
\, &AdS&Block I with $r>0$ &Block I with $r<0$ & Block II\\ \hline
$qe<0$ & + & - & + & - \\ \hline
$qe>0$ & - & + & - & + \\ \hline
\end{tabular}
\end{center}

\vspace{10pt}

Both of the table above assume a particle outside of $\Omega$. Inside of $\Omega$, we must negate every sign. Additionally, the $t'$ table assumes $r^2>a^2$. In particular, this will always be false in a region around $r=0$. Inside of this region, we must again negate the signs in the above $t'$ table. 
\subsection{Colatitude Coordinate Evolution}
We begin our examination of the colatitude coordinate evolution by noting a difference in how particles may approach the center of the singularity. Here we will have both $r=0$ and $\theta=\frac{\pi}{2}$. We begin by defining a new first integral, $\mathcal{Q}$, dependent upon the others, which appears to simply be a shifted form of Carter's constant: 
\begin{equation}
\mathcal{K}=\mathcal{Q}+\Xi^2(\mathcal{L}-a\mathcal{E})^2.
\end{equation}
We will refer to both $\mathcal{Q}$ and $\mathcal{K}$ as Carter's constant. Next, we examine our first integral equation $R(r)|_{r=0}$
\begin{align*}
R(0)&=\Delta_r(0)(-K)-\Xi^2\mathcal{P}(0)^2 \\
&=(a^2+e^2)(\mathcal{Q}+\Xi^2(\mathcal{L}-a\mathcal{E})^2)-\Xi^2(a\mathcal{L}-a^2\mathcal{E}^2)^2 \\
&=(a^2+e^2)\mathcal{Q}+e^2\Xi^2(\mathcal{L}-a\mathcal{E})^2.
\end{align*}
Enforcing the non-negativity of this expression only requires 
\begin{equation*}
\mathcal{Q}\geq \frac{-e^2\Xi^2(\mathcal{L}-a\mathcal{E})^2}{(a^2+e^2)}.
\end{equation*}
Next, we examine
\begin{align*}
\Theta\left(\frac{\pi}{2}\right)&=\Delta_\theta\left(\frac{\pi}{2}\right)\left(\mathcal{K}-qa^2C^2\left(\frac{\pi}{2}\right)\right)-\frac{\Xi^2\mathcal{D}^2}{S^2(\frac{\pi}{2})} \\
&=\mathcal{Q}+\Xi^2(\mathcal{L}-a\mathcal{E})^2-\Xi^2(\mathcal{L}-a\mathcal{E})^2 \\
&=\mathcal{Q},
\end{align*}
which shows we must only have $\mathcal{Q}\geq 0$. This is in stark contrast to the Kerr-Newman-de Sitter case, where we not only require $\mathcal{K}=0$, but we also require $\mathcal{L}=a\mathcal{E}$. Note that this result generalizes to the instanton containing $\theta$-horizons. 

Next, we aim to separate our $\mathcal{Q}$ function into a rotational kinetic energy function and a rotational potential energy function. To do this, we make use of the identity
\begin{equation*}
\Delta_\theta(\mathcal{L}-a\mathcal{E})^2-\frac{(\mathcal{L}-a\mathcal{E}S^2)^2}{S^2}=-a^2C^2\left(\frac{\mathcal{L}^2}{a^2S^2}+L(\mathcal{L}-a\mathcal{E})^2-\mathcal{E}^2\right),
\end{equation*}
which can be seen by expanding terms and using the definition of $\Delta_\theta$. It follows that 
\begin{align*}
\mathcal{Q}&=\mathcal{K}-\Xi^2(\mathcal{L}-a\mathcal{E})^2 \\
&=\mathfrak{q}a^2C^2+\frac{q}{\Delta_\theta}\left[\Sigma^2\theta'^2+\frac{\Xi^2\mathcal{D}^2}{S^2}\right]-\frac{\Xi^2\Delta_\theta}{\Delta_\theta}(\mathcal{L}-a\mathcal{E})^2 \\
&=\mathfrak{q}a^2C^2+\frac{1}{\Delta_\theta}\left[\Sigma^2\theta'^2+\Xi^2\left(\frac{(\mathcal{L}-a\mathcal{E}S^2)^2}{S^2}-\Delta_\theta(\mathcal{L}-a\mathcal{E})^2\right)\right] \\
&=\mathfrak{q}a^2C^2+\frac{1}{\Delta_\theta}\left[\Sigma^2\theta'^2+\Xi^2a^2C^2\left(L(\mathcal{L}-a\mathcal{E})^2-\mathcal{E}^2+\frac{\mathcal{L}}{a^2S^2}\right)\right]. 
\end{align*}
which we decompose into
\begin{equation}
\mathcal{Q}=T+V,
\end{equation}
where we define
\begin{align}
T((r,\theta),\theta'):=&\frac{\Sigma^2\theta'^2}{\Delta_\theta} \\
V(\theta):=&a^2C^2\left[\mathfrak{q}+\frac{\Xi^2}{\Delta_\theta}\left(L(\mathcal{L}-a\mathcal{E})^2+\frac{\mathcal{L}^2}{a^2S^2}-\mathcal{E}^2\right)\right],
\end{align}
as our \textit{rotational kinetic and potential energy functions} respectively. In the case of $\Lambda<0$, $L(\mathcal{L}-a\mathcal{E})^2-\mathcal{E}^2$ is a negative constant, which we will denote $\chi$.  In our $\theta$-horizonless instanton, we always have $\Delta_\theta>0$, making our kinetic energy function strictly positive. 
\subsubsection{Particles with Angular Momentum}
We first consider particles for which $\mathcal{L}\neq0$. Our potential energy function has an obvious root at $\theta=\frac{\pi}{2}$. It is clear that
\begin{equation*}
\lim_{\theta \to0 ,\pi}V(\theta)=+\infty.
\end{equation*}
However, due to the constancy of $\mathcal{Q}$, the potential energy may never truly reach infinite values. The potential energy function has an obvious root at $\theta=\frac{\pi}{2}$. Additional roots of this function will be solutions of the equation 
\begin{equation*}
S^2=\frac{\mathcal{L}^2}{a^2}\left(\frac{\Xi^2}{-\mathfrak{q}\Delta_\theta-\chi\Xi^2}\right).
\end{equation*}
A sufficient condition for the existence of roots thus becomes 
\begin{equation*}
0<\frac{\Xi^2}{-\mathfrak{q}\Delta_\theta-\chi\Xi^2}<\frac{a^2}{\mathcal{L}^2}.
\end{equation*}
When these two equations are satisfied, we will have two additional roots symmetric on either side of the equator, which we will label $\theta_{R-}$ and $\theta_{R+}$. Taking the derivative of the potential with respect to $\theta$ gives
\begin{equation*}
\frac{dV}{d\theta}=-2a^2CS\left[\mathfrak{q}+\left(\frac{\Xi^2}{\Delta_\theta^2}\left(\chi+\frac{\mathcal{L}^2}{a^2S^2}\right)+\frac{C^2\Xi^2\mathcal{L}^2}{\Delta_\theta a^2S^4}\right)\right],
\end{equation*}
which has limits
\begin{equation*}
\lim_{\theta\to0}=-\infty \hspace{1cm} \lim_{\theta\to\pi}=+\infty,
\end{equation*}
in accordance with our previous results. These results produce the same particle dynamics as in KNdS spacetimes or KNAdS spacetimes without $\theta$-horizons.
\begin{center}
\includegraphics[width=3.5in, height=3in]{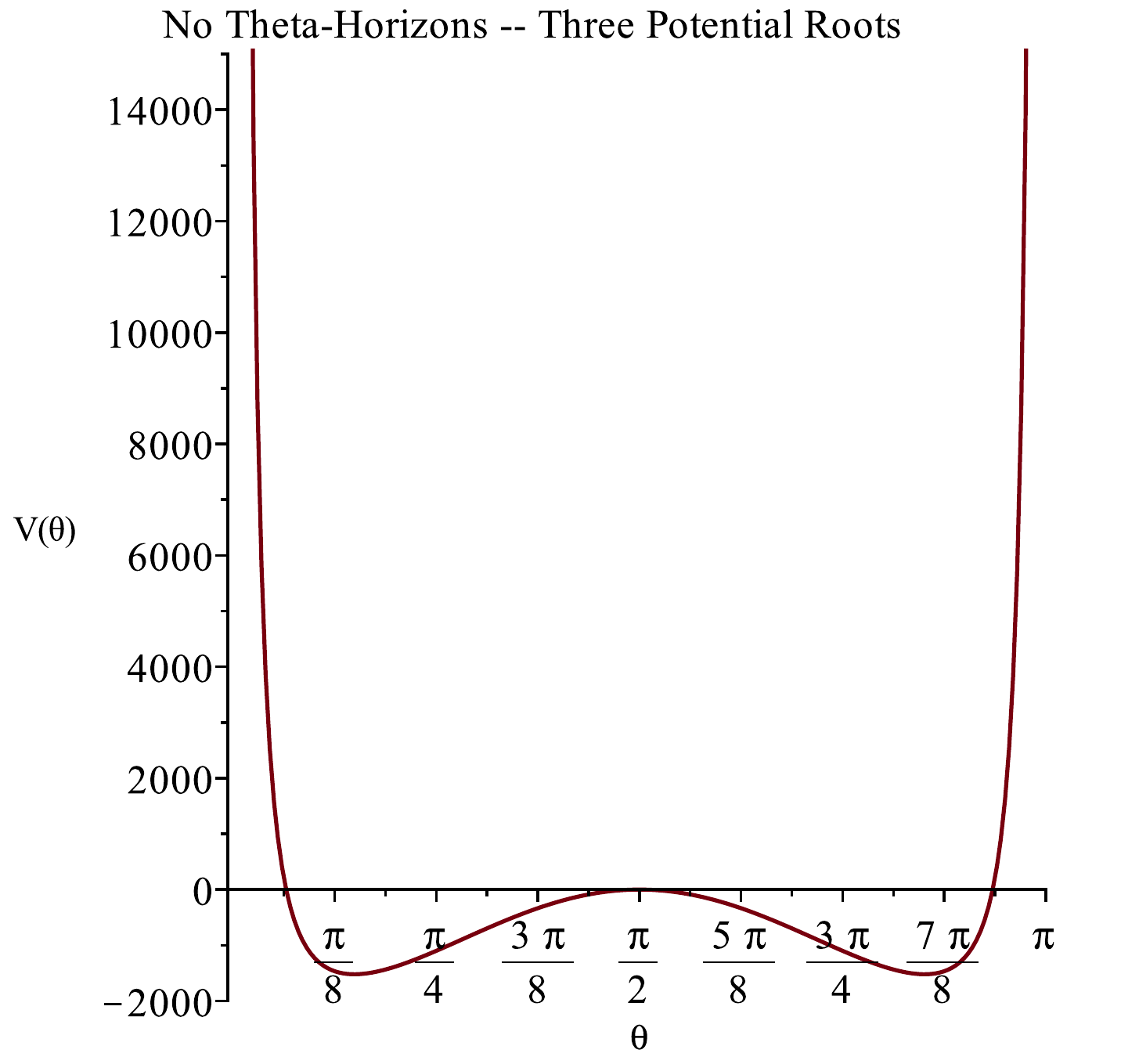}
\end{center}
Labelling the angular location of the additional minima $\theta_\pm^*$, we summarize the dynamics as follows:
\begin{itemize}
\item If we have additional roots
	\begin{itemize}
	\item $\mathcal{Q}>0 \Rightarrow \theta$ \textit{oscillates symmetrically about $\frac{\pi}{2}$}
\item $\mathcal{Q}=0 \Rightarrow \theta$, \textit{lies at an unstable equilibrium at $\frac{\pi}{2}$ or asymptotically approaches the equator from a negative potential}
\item $\mathcal{Q}_{min}<\mathcal{Q}<0 \Rightarrow \theta$ \textit{oscillates in one of the potential wells on either side of the equator}
\item $\mathcal{Q}=\mathcal{Q}_{min}\Rightarrow \theta$ \textit{lies in a stable equilibrium at $\theta_-^*$ or $\theta_+^*$}
	\end{itemize}
\item If we only have one root
	\begin{itemize}
	\item $\mathcal{Q}>0 \Rightarrow \theta$ \textit{oscillates symmetrically about $\frac{\pi}{2}$}
	\item $\mathcal{Q}=0 \Rightarrow \theta$ \textit{lies at a stable equilibrium at $\frac{\pi}{2}$}
	\end{itemize}
\end{itemize}
\subsubsection{Particles Without Angular Momentum}
Now we consider particles for which $\mathcal{L}=0$. These particles have the simpler potential energy function
\begin{equation}
V(\theta)=a^2C^2\left(\mathfrak{q}-\frac{\Xi^3\mathcal{E}^2}{\Delta_\theta}\right),
\end{equation}
with derivative
\begin{equation*}
\frac{dV}{d\theta}=-2a^2CS\left(\mathfrak{q}-\frac{\Xi^3\mathcal{E}^2}{\Delta_\theta^2}\right).
\end{equation*}
In contrast to particles with angular momentum, as these particles approach the poles, their potential ceases to diverge. Rather, it will approach the value 
\begin{equation*}
V(0,\pi)=a^2\left(\mathfrak{q}-\Xi^2\mathcal{E}^2\right).
\end{equation*}
This shows that if $\mathfrak{q}>\Xi^2\mathcal{E}^2$, we will have a positive potential at the poles, and the function will loosely resemble an upwards facing parabola. If $\mathfrak{q}$ is lightlike, timelike, or spacelike with $\mathfrak{q}<\Xi^2\mathcal{E}^2$, we will have negative potential at the poles, creating a shape resembling a downwards facing parabola. Clearly, we still have a potential root at $\theta=\frac{\pi}{2}$. Additional roots will be zeros of the equation
\begin{equation*}
\cos\theta=\pm\sqrt{\frac{1}{La^2}\left(1-\frac{\Xi^3\mathcal{E}^2}{\mathfrak{q}}\right)},
\end{equation*}
which in turn implies both
\begin{equation*}
\sqrt{\frac{\mathfrak{q}}{\Xi^3}}<\mathcal{E} \hspace{1cm} \left(1-\Xi\mathcal{E}\sqrt{\frac{\Xi}{\mathfrak{q}}}\right)\geq La^2.
\end{equation*}
\includegraphics[width=2.5in, height=2.5in]{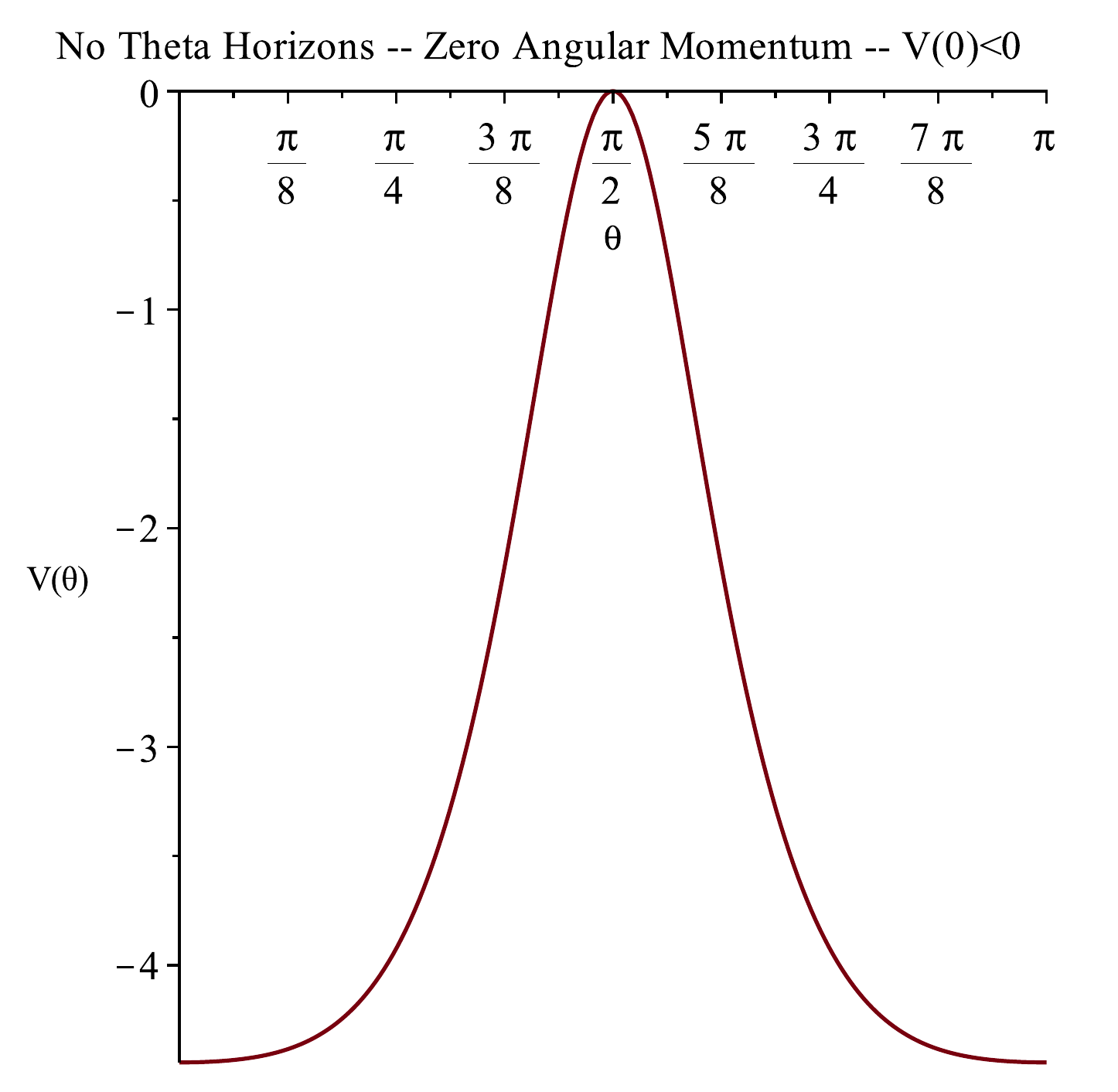}
\includegraphics[width=2.5in, height=2.5in]{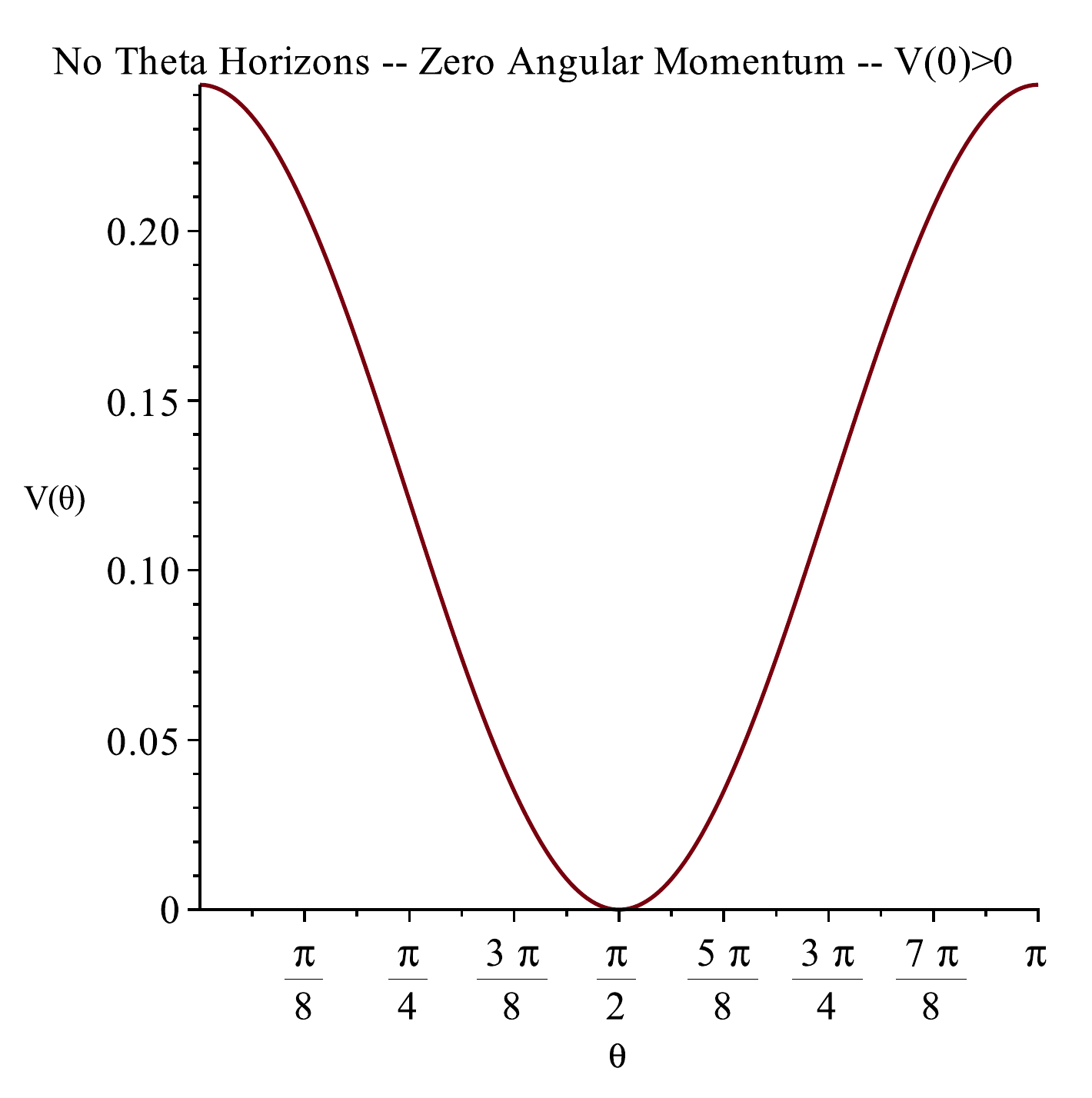}
This implies that only spacelike particles may have additional minima when $\mathcal{L}=0$. We summarize these dynamics as follows:
\begin{itemize}
\item If $V(0)>0$
	\begin{itemize}
	\item If one root
		\begin{itemize}
		\item $0<\mathcal{Q}<V_{max}\Rightarrow$ \textit{The behaviour will be the same as the case of particles with one potential root and angular momentum}
		\item $\mathcal{Q}\geq V_{max}\Rightarrow$ \textit{The particle will have a full range of motion between and on the poles}
		\end{itemize}
	\item If multiple roots
		\begin{itemize}
		\item $V_{min}\leq\mathcal{Q}<V_{max}\Rightarrow$ \textit{The behaviour will be the same as the case of particles with multiple potential roots and angular momentum}
		\item $V_{max}\leq \mathcal{Q}\Rightarrow$ \textit{The particle will have a full range of motion between and on the poles}
		\end{itemize}				
	\end{itemize}
\item If $V(0)<0$
	\begin{itemize}
	\item If one root
		\begin{itemize}
		\item $\mathcal{Q}\geq 0\Rightarrow$ \textit{The particle will have a full range of motion between and on the poles}
		\item $V_{min}\leq\mathcal{Q}<0\Rightarrow$ \textit{The particle will be able to move between $V=\mathcal{Q}$ and the poles}
		\end{itemize}
	\item If multiple roots
		\begin{itemize}
		\item $\mathcal{Q}\geq V_{max}\Rightarrow$ \textit{The particle will have a full range of motion between and on the poles}
		\item $0\leq\mathcal{Q}<V_{max}\Rightarrow$ \textit{The particle may move between $V=\mathcal{Q}$ and the poles, as well as the areas symmetric on either side of equator from $\theta=\frac{\pi}{2}$ to $V=\mathcal{Q}$}
		\item $V_{min}\leq\mathcal{Q}<0\Rightarrow$ \textit{The particle will be able to move between $V=\mathcal{Q}$ and the poles}
		\end{itemize}
	\end{itemize}		
\end{itemize}

\subsection{Radial Coordinate Evolution}
Finally, we briefly consider some aspects the radial evolution of a particle. We start with our radial equation 
\begin{equation*}
R(r)=\Sigma^2r'^2=\Delta_r(\mathfrak{q}r^2-\mathcal{K})-\Xi^2\mathcal{P}^2,
\end{equation*}
which is strictly non-negative. This implies we must have 
\begin{equation*}
(\mathfrak{q}r^2-\mathcal{K})\geq\Xi^2\mathcal{P}^2.
\end{equation*}
However, recall $\mathcal{K}=\mathfrak{q}r^2-\Sigma\langle\gamma'_\Pi,\gamma'_\Pi\rangle$. Substituting this into our inequality gives
\begin{equation}
\Sigma\langle\gamma'_\Pi,\gamma'_\Pi\rangle\geq\Xi^2\mathcal{P}^2,
\end{equation}
which is necessarily false whenever the quantity on the left hand side is negative. This quantity is \textit{always negative inside of Block I.} It appears, rather remarkably, that the entirety of the instanton between the event horizon and the inner event horizon is an inadmissible range of radii, regardless of whether you are inside of $\Omega$ or not.
\section{KNdS Instanton Particle Dynamics}
The Kerr-Newman-de Sitter Instanton creates a geometry much like that of the Kerr-Newman-anti-de Sitter spacetime, with both geometries containing $\theta$-horizons arising from the vanishing of the $\Delta_\theta$ term in the metric. This only occurs if the our angular momentum satisfies
\begin{equation}
a>a_{crit}=\sqrt{\frac{3}{\Lambda}}.
\end{equation}
If this angular momentum requirement is met, we shall call it a \textit{fast rotating instanton}, and if it is not, it shall be called a \textit{slowly rotating instanton}. Note that as $a\to\infty$, we have $\theta_\pm\to\frac{\pi}{2}$. Inside of these $\theta$-cones, as we shall call them, $\Delta_\theta$ is strictly negative. We begin by recalling the two expressions for Carter's constant:
\begin{align*}
\mathcal{K}&=\mathfrak{q}a^2C^2+\Sigma\langle\gamma'_\perp \gamma'_\perp\rangle \\
\mathcal{K}&=\mathfrak{q}r^2-\Sigma\langle\gamma'_\Pi\gamma'_\Pi\rangle.
\end{align*}
Using the $\gamma'_\perp$ formulation of $\mathcal{K}$, we consider the possible values $\mathcal{K}$ may take on, both inside and outside of the $\theta$-cones. In the slower rotating instanton lacking $\theta$-horizons, the possible values of $\mathcal{K}$ will be exactly those in Table 3.2. In the instanton with $\theta$-horizons, the possible values of $\mathcal{K}$ are summarized in tables 4.1 and 4.2:

\vspace{10pt}

\begin{center}
Table 4.1: \textit{Possible $\mathcal{K}$ values --- Fast Rotating Instanton --- Outside Cones}\vspace{5mm}
\begin{tabular}{c|c}
$\mathfrak{q}<0$ & $\mathcal{K}\geq \mathfrak{q}a\cos^2\theta_\pm$ with equality if a principal orbit at a $\theta$-horizon \\ \hline
$\mathfrak{q}=0$ & $\mathcal{K}=0$ if principal and $\mathcal{K}=\Sigma\langle\gamma'_\Pi,\gamma'_\Pi\rangle$ otherwise \\ \hline
$\mathfrak{q}>0$ & $\mathcal{K}\geq 0$ with equality if a principal equatorial orbit
\end{tabular}
\end{center} 

\vspace{10pt}

\begin{center}
Table 4.2: \textit{Possible $\mathcal{K}$ values --- Fast Rotating Instanton --- Inside Cones}\vspace{5mm}
\begin{tabular}{c|c}
$\mathfrak{q}<0$ & $\mathcal{K}< \mathfrak{q}a^2\cos^2\theta_\pm$ while approaching this value on the $\theta$-horizon \\
$\,$ & assumes a principal orbit \\ \hline
$\mathfrak{q}=0$ & $\mathcal{K}=0$ if principal and $\mathcal{K}=\Sigma\langle\gamma'_\Pi,\gamma'_\Pi\rangle$ otherwise \\ \hline
$\mathfrak{q}>0$ & $\mathcal{K}\leq \mathfrak{q}a^2$ with equality if a principal orbit on a pole 
\end{tabular}
\end{center}

\vspace{10pt}

In the preceding tables, the values of $\mathcal{K}$ for $\mathfrak{q}=0$ depend upon the following (Table 4.3):

\vspace{10pt}

\begin{center}
Table 4.3: \textit{Causal Character of $\langle\gamma'_\Pi,\gamma'_\Pi\rangle$ for KNdSI}\vspace{2mm}
\begin{tabular}{c|c}
$\langle\gamma'_\Pi,\gamma'_\Pi\rangle<0$ & dS$_\pm$ Block outside $\Omega$ -- Block I or Block III inside $\Omega$ \\ 
$\,$ &   Block II outside $\Omega$ \\ \hline
$\langle\gamma'_\Pi,\gamma'_\Pi\rangle>0$ & dS$_\pm$ Block inside $\Omega$ -- Block I or Block III outside $\Omega$ \\ 
$\,$ &   Block II inside $\Omega$ \\ \hline
\end{tabular}
\end{center}

\vspace{10pt}

\subsection{Time and Longitude Coordinate Evolution}
Again, we consider the case of lazy particles. In the slowly rotating instanton, the particles will move the same along the $t$ and $\varphi$ axes as in the KNAdSI instanton. However, in the fast rotating instanton, we have the new coordinate evolution as follows:

\begin{center}
Table 4.4: \textit{Direction of travel along the $\varphi'$ axis for $\mathcal{L}=\mathcal{E}=0$ when $\Lambda>0$ and $a>a_{crit}$} \\\vspace{2mm}
\begin{tabular}{c|c|c|c|c|c}
\, & dS$_\pm$ & Block I & Block II & Block III with $r>0$ & Block III with $r<0$ \\ \hline
$qe<0$ & $\mp$ & + & - & + & - \\ \hline
$qe>0$ & $\pm$ & - & + & - & + \\ \hline
\end{tabular}
\end{center}

\vspace{10pt}

\begin{center}
Table 4.5: \textit{Direction of travel along the $t'$ axis for $\mathcal{L}=\mathcal{E}=0$ when $\Lambda>0$ and $a>a_{crit}$} \\ \vspace{2mm}
\begin{tabular}{c|c|c|c|c|c}
\, & dS$_\pm$ &Block I & Block II & Block III with $r>0$ & Block III with $r<0$ \\ \hline
$qe<0$ & $\pm$ & - & + & - & + \\ \hline
$qe>0$ & $\mp$ & + & - & + & - \\ \hline
\end{tabular}
\end{center}

\vspace{10pt}

Again, we are assuming these particles lie outside of $\Omega$. If located inside of $\Omega$, one simply negates every sign. The assumption in the $t'$ table that $r^2>a^2$ is worse in this instanton, due to its large $a$. Everywhere $a^2>r^2$, we must again negate the signs in the $t'$ table.
\subsection{Colatitude Coordinate Evolution}
This subsection will focus on the evolution of the $\theta$ coordinate of particles in fast rotating instantons. In slowly rotating KNdSI instanton, the dynamics are nearly identical to KNAdSI, and are not worth exploring in detail. On the other hand, the existence of $\theta$-horizons in the fast rotating instanton, which act as impenetrable barriers of the $\theta$-coordinate (as you'd need infinite kinetic energy to cross one), create much more interesting behaviour. Recall respectively the rotational kinetic and potential energy functions
\begin{align*}
T((r,\theta),\theta')&=\frac{\Sigma^2\theta'^2}{\Delta_\theta} \\
V(\theta)&=a^2C^2\left[\mathfrak{q}+\frac{\Xi^2}{\Delta_\theta}\left(L(\mathcal{L}-a\mathcal{E})^2+\frac{\mathcal{L}}{a^2S^2}-\mathcal{E}^2\right)\right].
\end{align*}
\subsubsection{Particles With Angular Momentum}
First we will consider the case $\mathcal{L}\neq0$. With the $\theta$-horizons being zeroes of the $\Delta_\theta$ function, it will be important to understand how the function 
\begin{equation*}
\Psi:=L(\mathcal{L}-a\mathcal{E})^2+\frac{\mathcal{L}}{a^2S^2}-\mathcal{E}^2
\end{equation*}
behaves in order to understand behaviour of particles approaching the $\theta$-horizons. Setting $\Psi=0$ and rearranging gives
\begin{align*}
0&=L\mathcal{L}^2-2L\mathcal{L}a\mathcal{E}+La^2\mathcal{E}^2-\mathcal{E}^2+\frac{\mathcal{L}^2}{a^2S^2} \\
&=\mathcal{E}^2\Xi+2\mathcal{L}a\mathcal{E}-\mathcal{L}^2\left(L+\frac{1}{a^2S^2}\right).
\end{align*}
Solving this quadratic equation in $\mathcal{E}$ then gives 
\begin{equation*}
\mathcal{E}=\frac{-2L\mathcal{L}a\pm \sqrt{4L^2\mathcal{L}^2a^2+4\Xi\mathcal{L}^2\left(L+\frac{1}{a^2S^2}\right)}}{2\Xi}.
\end{equation*}
We can see this quadratic function in $\mathcal{E}$ creates a downwards facing parabola. In order for this parabola to have real roots, we require
\begin{align*}
0&\leq L^2a^2+\frac{\Xi(La^2S^2+1)}{a^2S^2} \\
&\leq \frac{La^2(S^2-1)+1}{a^2S^2} \\
&\leq \frac{\Delta_\theta}{a^2S^2}.
\end{align*}
This shows that $\Psi$ may only have negative values inside of the $\theta$-cones. Outside of the cones, where $\Delta_\theta>0$ we may have $\Psi$ be positive or negative. Analysis of $V(\theta)$ shows that we have the following limiting behaviour:

\vspace{10pt}

\begin{center}
Table 4.6: \textit{Limiting Behaviour of $V(\theta)$ for $\Lambda>0$ and $a>a_{crit}$}\vspace{5mm}
\begin{tabular}{c|c}
limit as $\theta$ approaches & $V(\theta)$ \\ \hline
0, $\pi$ & $-\infty$ \\ \hline
$\theta_\pm^\pm$ (Inside of Cones) & $+\infty$ \\ \hline
$\theta_\pm^\mp$ with $\Psi>0$ (Outside of Cones) & $+\infty$ \\ \hline
$\theta_\pm^\mp$ with $\Psi<0$ (Outside of Cones) & $-\infty$ \\ \hline
\end{tabular}
\end{center}

\vspace{10pt}

This clearly implies the existence of a mandatory potential root inside of the $\theta$-cones. Outside of the cones, additional roots will be solutions of the equation  
\begin{equation}
S^2=\frac{\mathcal{L}^2}{a^2}\left(\frac{\Xi^2}{\mathcal{E}^2\Xi^2-L(\mathcal{L}-a\mathcal{E})^2\Xi^2-\Delta_\theta\mathfrak{q}}\right).
\end{equation}
If the angular location of these solutions, $\theta_{R\pm}$, lies in the range 
\begin{equation*}
0<\theta_-<\theta_{R-}<\frac{\pi}{2}<\theta_{R+}<\theta_+<\pi,
\end{equation*}
then we will have 5 potential roots in our instanton. Note that inside of the cones, only $T<0$ is possible, implying a \textit{negative kinetic energy} here. Using the constancy of $\mathcal{Q}$, we analyse the dynamics arising from the various potential energy configurations that are possible. 
\begin{itemize}
\item Inside of the cones
	\begin{itemize}
	\item $\mathcal{Q}\geq0\Rightarrow$ \textit{the particle may approach the $\theta$-horizon asymptotically, with $V(\theta)\to\infty$ and $T\to - \infty$. Going towards the poles, the particle may only approach the $\theta$ value at which $V(\theta)=\mathcal{Q}$. At this point we will have $T=0\Rightarrow\theta'^2=0$.}
	\item $\mathcal{Q}<0\Rightarrow$ \textit{The particle may again approach the $\theta$-horizon asymptotically, provided that $T\to -\infty$ faster than $V\to\infty$ to keep $\mathcal{Q}$ constant. The particle may approach the pole, but we have $V,T<0$ here, so they will be bounded by $\mathcal{Q}$, preventing the particle from actually reaching the pole unless we have $\mathcal{Q}=-\infty$.}
	\end{itemize}
\item Outside of the cones: One Root
	\begin{itemize}
	\item If $\Psi>0$ at $\theta$-horizons
		\begin{itemize}
		\item $\mathcal{Q}>0\Rightarrow$ \textit{the particle oscillates symmetrically about the equator. }
	\item $\mathcal{Q}=0\Rightarrow$ \textit{the particle lies at a stable equilibrium at the equator. }
		\end{itemize}
	\item If $\Psi<0$ at $\theta$-horizons
		\begin{itemize}
		\item If $\mathcal{Q}\geq 0$ \textit{The particle will have a full range of motion between the $\theta$-horizons, increasing its kinetic energy as it moves closer to the horizons, so as to keep $\mathcal{Q}$ constant.}
		\item If $\mathcal{Q}\leq 0$ \textit{The particle will have a range of motion between the $\theta$-horizons and where $V(\theta)=\mathcal{Q}$}
		\end{itemize}
	\end{itemize}
\item Outside of the cones: Three Roots
	\begin{itemize}
	\item If $\Psi>0$ at $\theta$-horizons
		\begin{itemize}
		\item $\mathcal{Q}>0 \Rightarrow \theta$ \textit{oscillates symmetrically about $\frac{\pi}{2}$}
\item $\mathcal{Q}=0 \Rightarrow \theta$ \textit{lies at an unstable equilibrium at $\frac{\pi}{2}$ or asymptotically approaches the equator from a negative potential}
\item $\mathcal{Q}_{min}<\mathcal{Q}<0 \Rightarrow \theta$ \textit{oscillates in one of the potential wells on either side of the equator}
\item $\mathcal{Q}=\mathcal{Q}_{min}\Rightarrow \theta$ \textit{lies in a stable equilibrium at $\theta_-^*$ or $\theta_+^*$}
		\end{itemize}
	\item If $\Psi<0$ at $\theta$-horizons
		\begin{itemize}
		\item $\mathcal{Q}\geq V_{max}\Rightarrow$ \textit{The particle will have a full range of motion between the $\theta$-horizons}
		\item $0\leq\mathcal{Q}<V_{max}\Rightarrow$ \textit{The particle may move between $V=\mathcal{Q}$ and the $\theta$-horizons, as well as the areas symmetric on either side of equator from $\theta=\frac{\pi}{2}$ to $V=\mathcal{Q}$}
		\item $V_{min}\leq\mathcal{Q}<0\Rightarrow$ \textit{The particle will be able to move between $V=\mathcal{Q}$ and the $\theta$-horizons}
		\end{itemize}
	\end{itemize}
\end{itemize}
\subsubsection{Particles Without Angular Momentum}
We now turn our attention to the case of $\mathcal{L}=0$. Recall our simplified potential function and its derivative
\begin{align*}
V(\theta)&=a^2C^2\left(\mathfrak{q}-\frac{\Xi^3\mathcal{E}^2}{\Delta_\theta}\right) \\
\frac{dV}{d\theta}&=-2a^2SC\left(\mathfrak{q}-\frac{\Xi^3\mathcal{E}^2}{\Delta_\theta^2}\right).
\end{align*}
Here, we will again have a negative potential at the poles if $\mathfrak{q}$ is timelike, lightlike, or spacelike with $\mathfrak{q}<\Xi^2\mathcal{E}^2$. If $\mathfrak{q}>\Xi^2\mathcal{E}^2$ then we will have a positive potential at the poles. Again, a clear potential root is present at $\theta=\frac{\pi}{2}$. Additional roots will be zeroes of the equation 
\begin{equation*}
\cos\theta=\pm\sqrt{\frac{1}{La^2}\left(1-\frac{\Xi^3\mathcal{E}^2}{\mathfrak{q}}\right)}.
\end{equation*}
For these additional roots to lie outside of the cones, we require that 
\begin{equation*}
0<\frac{\Xi^3\mathcal{E}^2}{\mathfrak{q}},
\end{equation*}
which implies that $\mathfrak{q}$ must be timelike. It is clear that we have the following limits
\begin{equation*}
\lim_{\theta\to\theta_\pm^\pm}=-\infty \hspace{1cm} \lim_{\theta\to\theta_\pm^\mp}=+\infty.
\end{equation*}
\includegraphics[width=2.5in, height=2.5in]{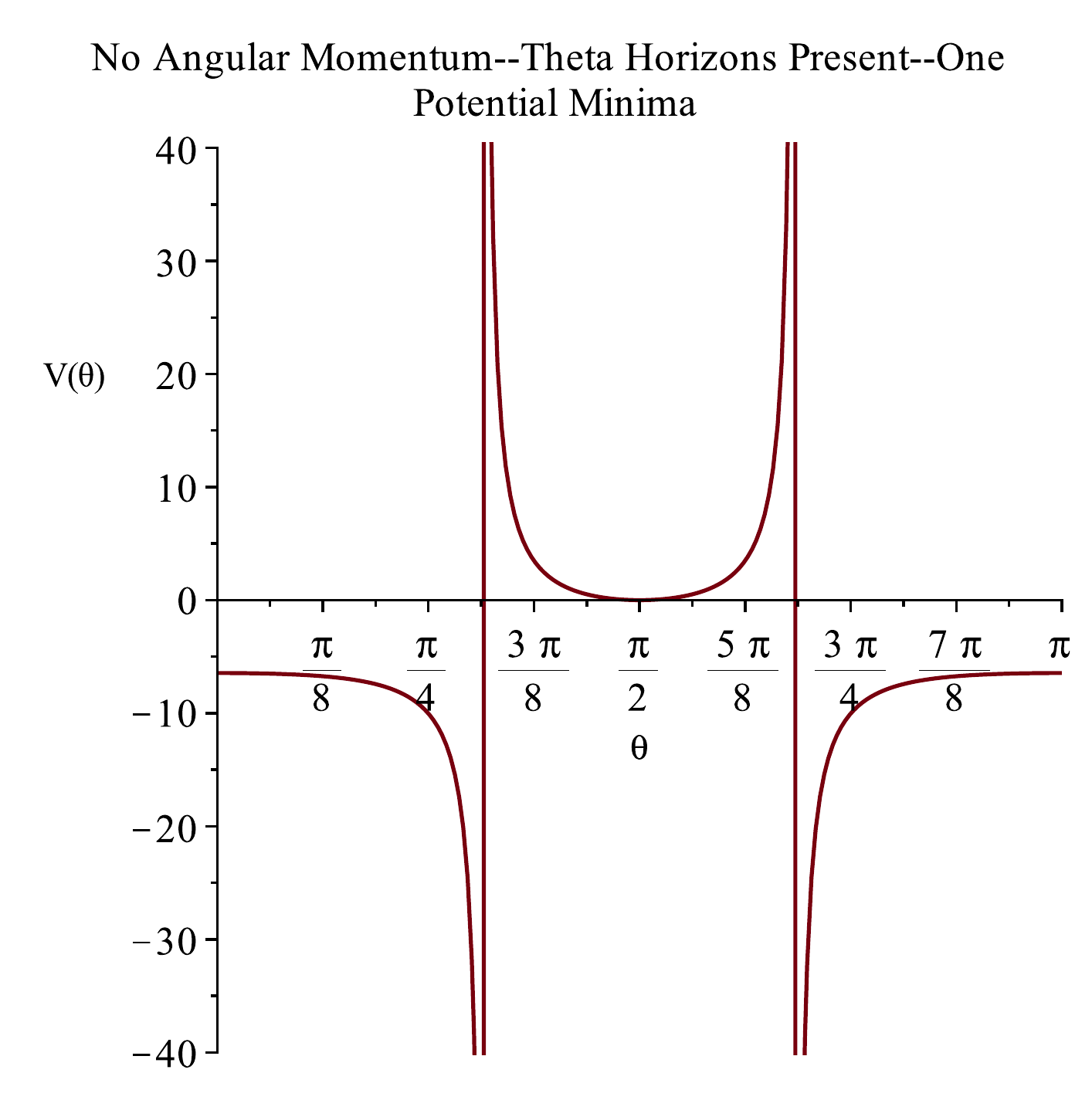}
\includegraphics[width=2.5in, height=2.5in]{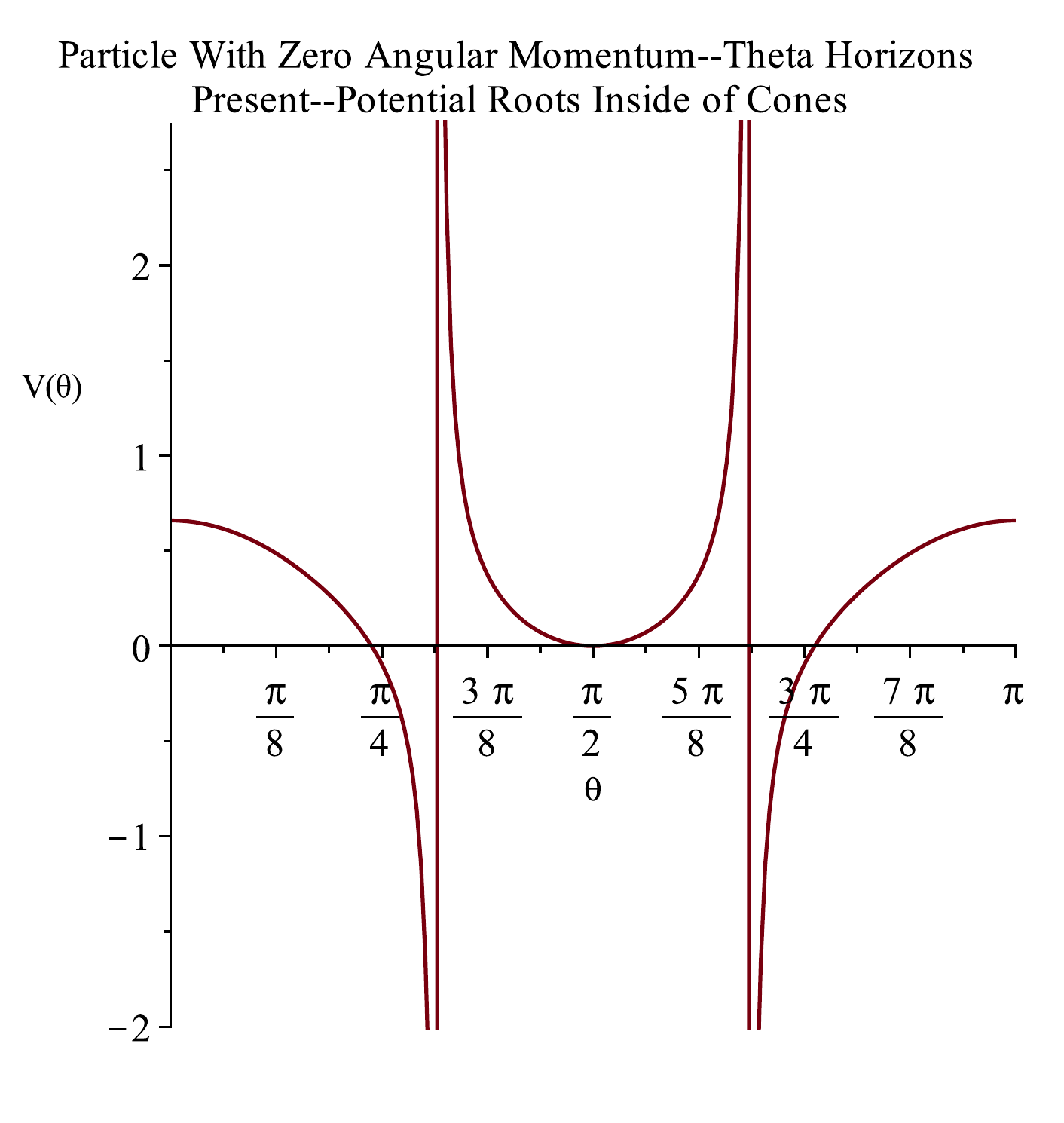}
We summarize these dynamical conclusions as follows
\begin{itemize}
\item Inside of the cones
	\begin{itemize}
	\item If $\mathfrak{q}<\Xi^2\mathcal{E}^2$
		\begin{itemize}
		\item $\mathcal{Q}=V_{max}\Rightarrow$\textit{ the particle starts at a pole and remains there}
		\item $\mathcal{Q}<V_{max}\Rightarrow$ \textit{The particle may reach the pole, and may approach the $\theta$-horizon until $\mathcal{Q}=V(\theta)$}
		\end{itemize}
	\item If $\mathfrak{q}>\Xi^2\mathcal{E}^2$
		\begin{itemize}
		 \item\textit{ We may have $\mathcal{Q}\leq V_{max}$ and may have the particle pass freely over the axis. The particle may only approach the $\theta$-horizons until $V(\theta)=\mathcal{Q}$}
		\end{itemize}
	\end{itemize}
\item Outside of the cones
	\begin{itemize}
	\item If only one root
		\begin{itemize}
		\item $\mathcal{Q}>0\Rightarrow$ \textit{the particle oscillates symmetrically about the equator. }
	\item $\mathcal{Q}=0\Rightarrow$ \textit{the particle lies at a stable equilibrium at the equator. }
		\end{itemize}
	\item If three roots
		\begin{itemize}
		\item $\mathcal{Q}>0 \Rightarrow \theta$ \textit{oscillates symmetrically about $\frac{\pi}{2}$}
\item $\mathcal{Q}=0 \Rightarrow \theta$, \textit{lies at an unstable equilibrium at $\frac{\pi}{2}$ or asymptotically approaches the equator from a negative potential}
\item $\mathcal{Q}_{min}<\mathcal{Q}<0 \Rightarrow \theta$ \textit{oscillates in one of the potential wells on either side of the equator}
\item $\mathcal{Q}=\mathcal{Q}_{min}\Rightarrow \theta$ \textit{lies in a stable equilibrium at $\theta_-^*$ or $\theta_+^*$}
		\end{itemize}			
	\end{itemize}
\end{itemize}

\subsection{Radial Coordinate Evolution}
Again, we briefly consider some aspects the radial evolution of a particle. We recall our inequality arising from the radial equation, and the substitution of Carter's constant: 
\begin{equation*}
\Sigma\langle\gamma'_\Pi,\gamma'_\Pi\rangle\geq\Xi^2\mathcal{P}^2.
\end{equation*}
The quantity on the left hand side of this equation is always negative for $r$ in the dS$_\pm$ blocks, or in Block II. It appears that these are inadmissible ranges of radii, regardless of whether or not the particle is inside of the singularity.

\vspace{10pt}

\noindent\textbf{Acknowledgements.}  We thank P. Chru\'sciel, M. Ghezelbash, N. Kamran, A. Shevyakov, J. Szmigielski, and Z. Stuchl\'{i}k for useful discussions related to this work.  The second named author acknowledges the support of a University of Saskatchewan New Faculty Recruitment Grant.  We thank the University of Saskatchewan and the University's College of Arts \& Science for support through their USRA program. 

\section{Appendix: Analysis of $\Delta_r$ Root Structure}

In the KN(A)dSI instanton, the radial horizons are determined by the roots of the quartic polynomial equation in $r$:

\begin{equation}
\Delta_r(r)=-Lr^4+r^2(La^2+1)-2Mr-(a^2+e^2),
\end{equation}

with $L=\frac{\Lambda}{3}$.

The possible root structures for a quartic parametric polynomial (that is, a polynomial with coefficients containing one or more parameters) are presented in \cite{Yang}. We obtain a system of three determinants whose signs will lend us information about the number of real roots in the instanton. These determinants take the form
\begin{center}
\begin{equation}
D1=8L(La^2+1)
\end{equation}
\begin{equation}
D2=16L^2(La^2+1)(-a^2-e^2)-72L^2M^2+4L(La^2+1)^3 
\end{equation}
\begin{multline}
D3=-256{L}^{3} \left( -{a}^{2}-{e}^{2} \right) ^{3}-128{L}^{2}
 \left( L{a}^{2}+1 \right) ^{2} \left( -{a}^{2}-{e}^{2} \right) ^{2} \\
+576{L}^{2} \left( L{a}^{2}+1 \right) {M}^{2} \left( -{a}^{2}-{e}^{2}
 \right) -432{L}^{2}{M}^{4} \\
 -16L \left( L{a}^{2}+1 \right) ^{4}
 \left( -{a}^{2}-{e}^{2} \right) +16L \left( L{a}^{2}+1 \right) ^{3}
{M}^{2}.
\end{multline}
\end{center}
For our radial horizon function to have two real roots and a pair of complex conjugate roots, we simply require that $D3<0$. Expanding $D3$, we obtain
\begin{multline*}
D3={L}^{5}{a}^{10}+{L}^{5}{a}^{8}{e}^{2}+{L}^{4}{M}^{2}{a}^{6}-4\,{L}^{4}
{a}^{8}-12\,{L}^{4}{a}^{6}{e}^{2}-8\,{L}^{4}{a}^{4}{e}^{4}-33\,{L}^{3}
{M}^{2}{a}^{4} \\
-36\,{L}^{3}{M}^{2}{a}^{2}{e}^{2}+6\,{L}^{3}{a}^{6}+22\,
{L}^{3}{a}^{4}{e}^{2}+32\,{L}^{3}{a}^{2}{e}^{4}+16\,{L}^{3}{e}^{6}-27
\,{L}^{2}{M}^{4}-33\,{L}^{2}{M}^{2}{a}^{2} \\
-36\,{L}^{2}{M}^{2}{e}^{2}-4
\,{L}^{2}{a}^{4}-12\,{L}^{2}{a}^{2}{e}^{2}-8\,{L}^{2}{e}^{4}+L{M}^{2}+
L{a}^{2}+L{e}^{2}.
\end{multline*}
Assuming $\Lambda<<1$, we discard all terms containing $\Lambda^n$ for $n\geq 2$. We are left with
\begin{equation*}
D3\approx L(a^2+M^2+e^2).
\end{equation*}
Thus, $D3<0$ for the $\Lambda<0$, and our KNAdSI instanton will only have two radial horizons. For the KNdSI instanton, we will have $D3>0$, and must continue investigating the determinants. It is clear that for KNdSI, we will have $D1>0$. \\
For four real roots, we require that $D2$ be strictly positive. Upon expanding $D2$, this inequality becomes
\begin{equation*}
{L}^{3}{a}^{6}-{L}^{2}{a}^{4}-4\,{L}^{2}{a}^{2}{e}^{2}-18\,L{M}^{2}-L{
a}^{2}-4\,L{e}^{2}+1>0.
\end{equation*}
This factors into the form 
\begin{equation*}
-18LM^2+(La^2-1)(L(La^2+re^2)-1)>0.
\end{equation*}
Making the approximation $(La^2-1)\approx -1$ we obtain
\begin{align*}
18M^2+La^2+4e^2&<\frac{1}{L} \\
\approx 18M^2+4e^2&<\frac{1}{L}.
\end{align*}
Since we have assumed $\Lambda << 1$, we take this to be a positive quantity. With this assumption, we will have four real roots in the KNdSI instanton.

Note that the existence of four real roots in the KNdSI instanton is incredibly sensitive to the magnitude of $\Lambda$. For systems with very large mass, the last approximation may not always be valid. In this case, there would exist no real roots whatsoever, and the black hole would cease to have any horizons at all.

Additionally, we observe that for both instantons, there will exist only one negative root of the radial horizon function. For the $\Lambda>0$ case, we put the radial horizon function into the form
\begin{equation}
Lr^4+(a^2+e^2)=r^2+La^2r^2-2Mr,
\end{equation}
which is strictly positive for $r<0$. Taking its derivative with respect to $r$, we obtain
\begin{equation*}
4L^3=2r+2La^2r-2M.
\end{equation*}
Both sides of this equation are strictly negative for $r<0$, hence both sides of (26) are strictly decreasing for $r<0$. For $r<<-1$ the left hand side of (26) is larger than the right hand side. However, at $r=0$, the left hand side is $(a^2+e^2)>0$, while the right hand side is $0$. Thus, both sides of this equation must cross at a unique point less than zero.

A similar analysis can be done for the $\Lambda<0$ case, again allowing only one negative root. Thus, there exists a region of Block III in the de Sitter instanton with $\Delta_r>0$, and a region with $\Delta_r<0$. Likewise in the anti-de Sitter instanton, we have a region of Block I with $\Delta_r>0$ and a region with $\Delta_r<0$.

\vfill
\pagebreak

\end{document}